\documentclass{amsproc}
\usepackage{amssymb}
\usepackage{amsmath}
\usepackage{eurosym}
\usepackage{amsfonts}

            \advance \topmargin by -\headheight \advance
            \topmargin by -\headsep
            \evensidemargin
            \oddsidemargin
\newtheorem{theorem}{Theorem}[section]

\newtheorem{corollary}[theorem]{Corollary}
\newtheorem{proposition}[theorem]{Proposition}

\numberwithin{equation}{section}
\input{tcilatex}

\begin{document}
\title[Summing Formulas For Generalized Tribonacci Numbers]{Summing Formulas
For Generalized Tribonacci Numbers}
\thanks{}
\author[Y\"{u}ksel~Soykan]{Y\"{u}ksel~Soykan}
\maketitle

\begin{center}
\textsl{Department of Mathematics, }

\textsl{Art and Science Faculty, }

\textsl{Zonguldak B\"{u}lent Ecevit University,}

\textsl{67100, Zonguldak, Turkey }

\textsl{e-mail: \ yuksel\_soykan@hotmail.com}
\end{center}

\textbf{Abstract. }In this paper, closed forms of the summation formulas for
generalized Tribonacci numbers are presented. Then, some previous results
are recovered as particular cases of the present results. As special cases,
we give summation formulas of Tribonacci, Tribonacci-Lucas, Padovan, Perrin,
Narayana and some other third order linear recurrance sequences. All the
summing formulas of well known recurrence sequences which we deal with are
linear except the cases Pell-Padovan and Padovan-Perrin.



\section{Introduction}

In this work, we investigate linear summation formulas of generalized
Tribonacci numbers. Some summing formulas of the Pell and Pell-Lucas numbers
are well known and given in [\ref{koshyfibbook2001}, \ref{koshypellbook2014}%
], see also [\ref{gokbassum2017}]. For linear sums of Fibonacci, Tribonacci,
Tetranacci, Pentanacci and Hexanacci numbers, see [\ref{hansen1978},\ref%
{soykan2019sumgenfibgagenfib}], [\ref{bib:frontczak2018},\ref{parpartez2011}%
], [\ref{soykan2018triposmatrix}, \ref{waddill1992}], [\ref%
{soykan2019penlinsum}], and [\ref{soykan2019hexanasum}] respectively.

First, in this section, we present some background about generalized
Tribonacci numbers. The generalized Tribonacci sequence $%
\{W_{n}(W_{0},W_{1},W_{2};r,s,t)\}_{n\geq 0}$ (or shortly $\{W_{n}\}_{n\geq
0}$) is defined as follows: 
\begin{equation}
W_{n}=rW_{n-1}+sW_{n-2}+tW_{n-3},\text{ \ \ \ \ }W_{0}=a,W_{1}=b,W_{2}=c,%
\text{ \ }n\geq 3  \label{equation:nbvcftyopuybnaeuo}
\end{equation}%
where $W_{0},W_{1},W_{2}\ $are arbitrary complex numbers and $r,s,t$ are
real numbers. The generalized Tribonacci sequence has been studied by many
authors, see for example [\ref{bib:bruce1984},\ref{catalani2002},\ref%
{bib:choi2013},\ref{bib:elia2001},\ref{bib:lin1988},\ref{bib:pethe1988},\ref%
{bib:scott1977},\ref{bib:shannon1972},\ref{bib:shannon1977},\ref%
{bib:spickerman1981},\ref{bib:yalavigi1971},\ref{bib:yalavigi1972},\ref%
{bib:yilmaz 2014},\ref{bib:waddill1991}].

The sequence $\{W_{n}\}_{n\geq 0}$ can be extended to negative subscripts by
defining%
\begin{equation*}
W_{-n}=-\frac{s}{t}W_{-(n-1)}-\frac{r}{t}W_{-(n-2)}+\frac{1}{t}W_{-(n-3)}
\end{equation*}%
for $n=1,2,3,...$ when $t\neq 0.$ Therefore, recurrence (\ref%
{equation:nbvcftyopuybnaeuo}) holds for all integer $n.$

If we set $r=s=t=1$ and $W_{0}=0,W_{1}=1,W_{2}=1$ then $\{W_{n}\}$ is the
well-known Tribonacci sequence and if we set $r=s=t=1$ and $%
W_{0}=3,W_{1}=1,W_{2}=3$ then $\{W_{n}\}$ is the well-known Tribonacci-Lucas
sequence.

In fact, the generalized Tribonacci sequence is the generalization of the
well-known sequences like Tribonacci, Tribonacci-Lucas, Padovan
(Cordonnier), Perrin, Padovan-Perrin, Narayana, third order Jacobsthal and
third order Jacobsthal-Lucas.\ In literature, for example, the following
names and notations (see Table 1) are used for the special case of $r,s,t$
and initial values.

\ \ \ \ \ 

Table 1 A few special case of generalized Tribonacci sequences.

\begin{tabular}{ccccc}
\hline
$\text{Sequences (Numbers)}$ &  & $\text{Notation}$ &  & OEIS [\ref%
{bib:sloane}] \\ \hline
$\text{Tribonacci}$ &  & $\{T_{n}\}=\{W_{n}(0,1,1;1,1,1)\}$ &  & A000073,
A057597 \\ 
$\text{Tribonacci-Lucas}$ &  & $\{K_{n}\}=\{W_{n}(3,1,3;1,1,1)\}$ &  & 
A001644, A073145 \\ 
$\text{third order }$Pell &  & $\{P_{n}^{(3)}\}=\{W_{n}(0,1,2;2,1,1)\}$ &  & 
A077939, A077978 \\ 
$\text{third order }$Pell-Lucas &  & $\{Q_{n}^{(3)}\}=\{W_{n}(3,2,6;2,1,1)\}$
&  & A276225, A276228 \\ 
$\text{third order modified }$Pell &  & $\{E_{n}^{(3)}\}=%
\{W_{n}(0,1,1;2,1,1)\}$ &  & A077997, A078049 \\ 
$\text{Padovan (Cordonnier)}$ &  & $\{P_{n}\}=\{W_{n}(1,1,1;0,1,1)\}$ &  & 
A000931 \\ 
$\text{Perrin (Padovan-Lucas)}$ &  & $\{E_{n}\}=\{W_{n}(3,0,2;0,1,1)\}$ &  & 
A001608, A078712 \\ 
$\text{Padovan-Perrin}$ &  & $\{S_{n}\}=\{W_{n}(0,0,1;0,1,1)\}$ &  & 
A000931, A176971 \\ 
$\text{Pell-Padovan}$ &  & $\{R_{n}\}=\{W_{n}(1,1,1;0,2,1)\}$ &  & A066983,
A128587 \\ 
$\text{Pell-Perrin}$ &  & $\{C_{n}\}=\{W_{n}(3,0,2;0,2,1)\}$ &  & - \\ 
$\text{Jacobsthal-Padovan}$ &  & $\{Q_{n}\}=\{W_{n}(1,1,1;0,1,2)\}$ &  & 
A159284 \\ 
$\text{Jacobsthal-Perrin (-Lucas)}$ &  & $\{D_{n}\}=\{W_{n}(3,0,2;0,1,2)\}$
&  & A072328 \\ 
$\text{Narayana}$ &  & $\{N_{n}\}=\{W_{n}(0,1,1;1,0,1)\}$ &  & A078012 \\ 
$\text{third order Jacobsthal}$ &  & $\{J_{n}^{(3)}\}=\{W_{n}(0,1,1;1,1,2)\}$
&  & A077947 \\ 
$\text{third order Jacobsthal-Lucas}$ &  & $\{j_{n}^{(3)}\}=%
\{W_{n}(2,1,5;1,1,2)\}$ &  & A226308 \\ \hline
\end{tabular}

Note that the sequence $\{C_{n}\}$ is't in the database of http://oeis.org [%
\ref{bib:sloane}], yet.

\section{Sum formulas of Generalized Tribonacci Numbers with Positive
Subscripts}

The following Theorem presents some linear summing formulas of generalized
Tribonacci numbers with positive subscripts.

\begin{theorem}
\label{theorem:firstbncvg}For $n\geq 0,$ we have the following formulas:

\begin{description}
\item[(a)] (Sum of the generalized Tribonacci numbers) If $r+s+t-1\neq 0,$
then 
\begin{equation*}
\sum_{k=0}^{n}W_{k}=\frac{\allowbreak
W_{n+3}+(1-r)W_{n+2}+(1-r-s)W_{n+1}-W_{2}+(r-1)W_{1}+(r+s-1)W_{0}}{r+s+t-1}.
\end{equation*}

\item[(b)] If $2s+2rt+r^{2}-s^{2}+t^{2}-1=\left( r+s+t-1\right) \left(
r-s+t+1\right) \neq 0$ then%
\begin{equation*}
\sum_{k=0}^{n}W_{2k}=\frac{%
\begin{array}{c}
(-s+1)W_{2n+2}+(t+rs)W_{2n+1}+(t^{2}+rt)W_{2n}+(-1+s)W_{2} \\ 
+(-t-rs)W_{1}+(-1+r^{2}-s^{2}+rt+2s)W_{0}%
\end{array}%
}{\left( r+s+t-1\right) \left( r-s+t+1\right) }
\end{equation*}%
and 
\begin{equation*}
\sum_{k=0}^{n}W_{2k+1}=\frac{%
\begin{array}{c}
(r+t)W_{2n+2}+(s-s^{2}+t^{2}+rt)W_{2n+1}+(t-st)W_{2n}+(-r-t)W_{2} \\ 
+(-1+s+r^{2}+rt)W_{1}+(-t+st)W_{0}%
\end{array}%
}{\left( r-s+t+1\right) \left( r+s+t-1\right) }.
\end{equation*}

\item[(c)] If $r+t\neq 0\wedge s=1$ then%
\begin{equation*}
\sum_{k=0}^{n}W_{2k}=\frac{1}{r+t}\left( W_{2n+1}+tW_{2n}-W_{1}+rW_{0}\right)
\end{equation*}%
and%
\begin{equation*}
\sum_{k=0}^{n}W_{2k+1}=\frac{1}{r+t}\left(
W_{2n+2}+tW_{2n+1}-W_{2}+rW_{1}\right) .
\end{equation*}%
Note that (c) is a special case of (b).
\end{description}
\end{theorem}

\textit{Proof.} \ 

\begin{description}
\item[(a)] Using the recurrence relation%
\begin{equation*}
W_{n}=rW_{n-1}+sW_{n-2}+tW_{n-3}
\end{equation*}%
i.e.%
\begin{equation*}
tW_{n-3}=W_{n}-rW_{n-1}-sW_{n-2}
\end{equation*}%
we obtain%
\begin{eqnarray*}
tW_{0} &=&W_{3}-rW_{2}-sW_{1} \\
tW_{1} &=&W_{4}-rW_{3}-sW_{2} \\
tW_{2} &=&W_{5}-rW_{4}-sW_{3} \\
&&\vdots \\
tW_{n-1} &=&W_{n+2}-rW_{n+1}-sW_{n} \\
tW_{n} &=&W_{n+3}-rW_{n+2}-sW_{n+1}.
\end{eqnarray*}%
If we add the equations by side by, we get%
\begin{equation*}
\sum_{k=0}^{n}W_{k}=\frac{\allowbreak
W_{n+3}+(1-r)W_{n+2}+(1-r-s)W_{n+1}-W_{2}+(r-1)W_{1}+(r+s-1)W_{0}}{r+s+t-1}.
\end{equation*}

\item[(b) and (c)] Using the recurrence relation%
\begin{equation*}
W_{n}=rW_{n-1}+sW_{n-2}+tW_{n-3}
\end{equation*}%
i.e.%
\begin{equation*}
rW_{n-1}=W_{n}-sW_{n-2}-tW_{n-3}
\end{equation*}%
we obtain%
\begin{eqnarray*}
rW_{3} &=&W_{4}-sW_{2}-tW_{1} \\
rW_{5} &=&W_{6}-sW_{4}-tW_{3} \\
&&\vdots \\
rW_{2n+1} &=&W_{2n+2}-sW_{2n}-tW_{2n-1}. \\
rW_{2n+3} &=&W_{2n+4}-sW_{2n+2}-tW_{2n+1}
\end{eqnarray*}%
Now, if we add the above equations by side by, we get%
\begin{equation}
r(-W_{1}+\sum_{k=0}^{n}W_{2k+1})=(W_{2n+2}-W_{2}-W_{0}+%
\sum_{k=0}^{n}W_{2k})-s(-W_{0}+\sum_{k=0}^{n}W_{2k})-t(-W_{2n+1}+%
\sum_{k=0}^{n}W_{2k+1}).  \label{equat:rtdynbcposd}
\end{equation}%
Similarly, using the recurrence relation%
\begin{equation*}
W_{n}=rW_{n-1}+sW_{n-2}+tW_{n-3}
\end{equation*}%
i.e.%
\begin{equation*}
rW_{n-1}=W_{n}-sW_{n-2}-tW_{n-3}
\end{equation*}%
we write the following obvious equations;%
\begin{eqnarray*}
rW_{2} &=&W_{3}-sW_{1}-tW_{0} \\
rW_{4} &=&W_{5}-sW_{3}-tW_{2} \\
rW_{6} &=&W_{7}-sW_{5}-tW_{4} \\
&&\vdots \\
rW_{2n} &=&W_{2n+1}-sW_{2n-1}-tW_{2n-2} \\
rW_{2n+2} &=&W_{2n+3}-sW_{2n+1}-tW_{2n}.
\end{eqnarray*}%
Now, if we add the above equations by side by, we obtain%
\begin{equation}
r(-W_{0}+\sum_{k=0}^{n}W_{2k})=(-W_{1}+\sum_{k=0}^{n}W_{2k+1})-s(-W_{2n+1}+%
\sum_{k=0}^{n}W_{2k+1})-t(-W_{2n}+\sum_{k=0}^{n}W_{2k}).
\label{equat:korksbnxcrtab}
\end{equation}%
Then, solving the system (\ref{equat:rtdynbcposd})-(\ref{equat:korksbnxcrtab}%
), the required results of (b) and (c) follow. 
\endproof%
\end{description}

Taking $r=s=t=1$ in Theorem \ref{theorem:firstbncvg} (a) and (b) (or (c)),
we obtain the following Proposition.

\begin{proposition}
If $r=s=t=1$ then for $n\geq 0$ we have the following formulas:

\begin{description}
\item[(a)] $\sum_{k=0}^{n}W_{k}=\frac{1}{2}(W_{n+3}-W_{n+1}-W_{2}+W_{0}).$

\item[(b)] $\sum_{k=0}^{n}W_{2k}=\frac{1}{2}(W_{2n+1}+W_{2n}-W_{1}+W_{0}).$

\item[(c)] $\sum_{k=0}^{n}W_{2k+1}=\frac{1}{2}%
(W_{2n+2}+W_{2n+1}-W_{2}+W_{1}).$
\end{description}
\end{proposition}

From the above Proposition, we have the following Corollary which gives
linear sum formulas of Tribonacci numbers (take $W_{n}=T_{n}$ with $%
T_{0}=0,T_{1}=1,T_{2}=1$).

\begin{corollary}
For $n\geq 0,$ Tribonacci numbers have the following properties.

\begin{description}
\item[(a)] $\sum_{k=0}^{n}T_{k}=\frac{1}{2}(T_{n+3}-T_{n+1}-1).$

\item[(b)] $\sum_{k=0}^{n}T_{2k}=\frac{1}{2}(T_{2n+1}+T_{2n}-1).$

\item[(c)] $\sum_{k=0}^{n}T_{2k+1}=\frac{1}{2}(T_{2n+2}+T_{2n+1}).$
\end{description}
\end{corollary}

Taking $W_{n}=K_{n}$ with $K_{0}=3,K_{1}=1,K_{2}=3$ in the above
Proposition, we have the following Corollary which presents linear sum
formulas of Tribonacci-Lucas numbers.

\begin{corollary}
For $n\geq 0,$ Tribonacci-Lucas numbers have the following properties.

\begin{description}
\item[(a)] $\sum_{k=0}^{n}K_{k}=\frac{1}{2}(K_{n+3}-K_{n+1}).$

\item[(b)] $\sum_{k=0}^{n}K_{2k}=\frac{1}{2}(K_{2n+1}+K_{2n}+2).$

\item[(c)] $\sum_{k=0}^{n}K_{2k+1}=\frac{1}{2}(K_{2n+2}+K_{2n+1}-2).$
\end{description}
\end{corollary}

Taking $r=2,s=1,t=1$ in Theorem \ref{theorem:firstbncvg} (a) and (b) (or
(c)), we obtain the following Proposition.

\begin{proposition}
If $r=2,s=1,t=1$ then for $n\geq 0$ we have the following formulas:

\begin{description}
\item[(a)] $\sum_{k=0}^{n}W_{k}=\frac{1}{3}\left( \allowbreak
W_{n+3}-W_{n+2}-2W_{n+1}-W_{2}+W_{1}+2W_{0}\right) .$

\item[(b)] $\sum_{k=0}^{n}W_{2k}=\frac{1}{3}\left(
W_{2n+1}+W_{2n}-W_{1}+2W_{0}\right) .$

\item[(c)] $\sum_{k=0}^{n}W_{2k+1}=\frac{1}{3}\left(
W_{2n+2}+W_{2n+1}-W_{2}+2W_{1}\right) .$
\end{description}
\end{proposition}

From the last Proposition, we have the following Corollary which gives
linear sum formulas of third-order Pell numbers (take $W_{n}=P_{n}^{(3)}$
with $P_{0}^{(3)}=0,P_{1}^{(3)}=1,P_{2}^{(3)}=2$).

\begin{corollary}
For $n\geq 0,$ third-order Pell numbers have the following properties:

\begin{description}
\item[(a)] $\sum_{k=0}^{n}P_{k}^{(3)}=\frac{1}{3}\allowbreak
(P_{n+3}^{(3)}-P_{n+2}^{(3)}-2P_{n+1}^{(3)}-1)$.

\item[(b)] $\sum_{k=0}^{n}P_{2k}^{(3)}=\frac{1}{3}%
(P_{2n+1}^{(3)}+P_{2n}^{(3)}-1).$

\item[(c)] $\sum_{k=0}^{n}P_{2k+1}^{(3)}=\frac{1}{3}%
(P_{2n+2}^{(3)}+P_{2n+1}^{(3)}).$
\end{description}
\end{corollary}

Taking $W_{n}=Q_{n}^{(3)}$ with $Q_{0}^{(3)}=3,Q_{1}^{(3)}=2,Q_{2}^{(3)}=6$
in the last Proposition, we have the following Corollary which presents
linear sum formulas of third-order Pell-Lucas numbers.

\begin{corollary}
For $n\geq 0,$ third-order Pell-Lucas numbers have the following properties:

\begin{description}
\item[(a)] $\sum_{k=0}^{n}Q_{k}^{(3)}=\frac{1}{3}\allowbreak
(Q_{n+3}^{(3)}-Q_{n+2}^{(3)}-2Q_{n+1}^{(3)}+2).$

\item[(b)] $\sum_{k=0}^{n}Q_{2k}^{(3)}=\frac{1}{3}%
(Q_{2n+1}^{(3)}+Q_{2n}^{(3)}+4).$

\item[(c)] $\sum_{k=0}^{n}Q_{2k+1}^{(3)}=\frac{1}{3}%
(Q_{2n+2}^{(3)}+Q_{2n+1}^{(3)}-2).$
\end{description}
\end{corollary}

From the last Proposition, we have the following Corollary which gives
linear sum formulas of third-order modified Pell numbers (take $%
W_{n}=E_{n}^{(3)}$ with $E_{0}^{(3)}=0,E_{1}^{(3)}=1,E_{2}^{(3)}=1$).

\begin{corollary}
For $n\geq 0,$ third-order modified Pell numbers have the following
properties:

\begin{description}
\item[(a)] $\sum_{k=0}^{n}E_{k}^{(3)}=\frac{1}{3}(\allowbreak
E_{n+3}^{(3)}-E_{n+2}^{(3)}-2E_{n+1}^{(3)})$

\item[(b)] $\sum_{k=0}^{n}E_{2k}^{(3)}=\frac{1}{3}%
(E_{2n+1}^{(3)}+E_{2n}^{(3)}-1)$

\item[(c)] $\sum_{k=0}^{n}E_{2k+1}^{(3)}=\frac{1}{3}%
(E_{2n+2}^{(3)}+E_{2n+1}^{(3)}+1).$
\end{description}
\end{corollary}

Taking $r=0,s=1,t=1$ in Theorem \ref{theorem:firstbncvg} (a) and (b) (or
(c)), we obtain the following Proposition.

\begin{proposition}
If $r=0,s=1,t=1$ then for $n\geq 0$ we have the following formulas:

\begin{description}
\item[(a)] $\sum_{k=0}^{n}W_{k}=W_{n+3}+W_{n+2}-W_{2}-W_{1}.$

\item[(b)] $\sum_{k=0}^{n}W_{2k}=W_{2n+1}+W_{2n}-W_{1}.$

\item[(c)] $\sum_{k=0}^{n}W_{2k+1}=W_{2n+2}+W_{2n+1}-W_{2}.$
\end{description}
\end{proposition}

From the last Proposition, we have the following Corollary which gives
linear sum formulas of Padovan numbers (take $W_{n}=P_{n}$ with $%
P_{0}=1,P=1,P_{2}=1$).

\begin{corollary}
For $n\geq 0,$ Padovan numbers have the following properties.

\begin{description}
\item[(a)] $\sum_{k=0}^{n}P_{k}=P_{n+3}+P_{n+2}-2.$

\item[(b)] $\sum_{k=0}^{n}P_{2k}=P_{2n+1}+P_{2n}-1.$

\item[(c)] $\sum_{k=0}^{n}P_{2k+1}=P_{2n+2}+P_{2n+1}-1.$
\end{description}
\end{corollary}

Taking $W_{n}=E_{n}$ with $E_{0}=3,E_{2}=0,E_{2}=2$ in the last Proposition,
we have the following Corollary which presents linear sum formulas of Perrin
numbers.

\begin{corollary}
For $n\geq 0,$ Perrin numbers have the following properties.

\begin{description}
\item[(a)] $\sum_{k=0}^{n}E_{k}=E_{n+3}+E_{n+2}-2.$

\item[(b)] $\sum_{k=0}^{n}E_{2k}=E_{2n+1}+E_{2n}.$

\item[(c)] $\sum_{k=0}^{n}E_{2k+1}=E_{2n+2}+E_{2n+1}-2.$
\end{description}
\end{corollary}

Taking $W_{n}=S_{n}$ with $S_{0}=0,S_{2}=0,S_{2}=1$ in the last Proposition,
we have the following Corollary which gives linear sum formulas of
Padovan-Perrin numbers.

\begin{corollary}
For $n\geq 0,$ Padovan-Perrin numbers have the following properties.

\begin{description}
\item[(a)] $\sum_{k=0}^{n}S_{k}=S_{n+3}+S_{n+2}-1.$

\item[(b)] $\sum_{k=0}^{n}S_{2k}=S_{2n+1}+S_{2n}.$

\item[(c)] $\sum_{k=0}^{n}S_{2k+1}=S_{2n+2}+S_{2n+1}-1.$
\end{description}
\end{corollary}

If $r=0,s=2,t=1$ then $\left( r-s+t+1\right) =\allowbreak 0$ so we can't use
Theorem \ref{theorem:firstbncvg} (b). In other words, the method of the
proof Theorem \ref{theorem:firstbncvg} (b) can't be used to find $%
\sum_{k=0}^{n}W_{2k}$ and $\sum_{k=0}^{n}W_{2k+1}.$ Therefore we need
another method to find them which is given in the following Theorem.

\begin{theorem}
\label{theorem:bayrambhospd}If $r=0,s=2,t=1$ then for $n\geq 0$ we have the
following formulas:

\begin{description}
\item[(a)] $\sum_{k=0}^{n}W_{k}=\frac{1}{2}\left(
W_{n+3}+W_{n+2}-W_{n+1}-W_{2}-W_{1}+W_{0}\right) .$

\item[(b)] $\sum_{k=0}^{n}W_{2k}=W_{2n+1}+\left( W_{2}-W_{1}-W_{0}\right)
n+W_{0}-W_{1}.$

\item[(c)] $\sum_{k=0}^{n}W_{2k+1}=\frac{1}{2}\left(
W_{2n+3}+W_{2n+2}-W_{2n+1}+2n\left( -W_{2}+W_{1}+W_{0}\right)
-W_{2}+W_{1}-W_{0}\right) .$
\end{description}
\end{theorem}

Proof.

\begin{description}
\item[(a)] Taking $r=0,s=2,t=1$ in Theorem \ref{theorem:firstbncvg} (a) we
obtain (a).

\item[(b) and (c)] Using the recurrence relation%
\begin{equation*}
W_{n}=2W_{n-2}+W_{n-3}
\end{equation*}%
we obtain%
\begin{eqnarray*}
\sum_{k=0}^{0}W_{2k} &=&W_{0} \\
\sum_{k=0}^{1}W_{2k} &=&W_{0}+W_{2}=W_{3}+W_{2}-2W_{1} \\
\sum_{k=0}^{2}W_{2k} &=&W_{0}+W_{2}+W_{4}=W_{5}+2W_{2}-3W_{1}-W_{0} \\
&&\vdots \\
\sum_{k=0}^{n}W_{2k} &=&W_{2n+1}+\left( W_{2}-W_{1}-W_{0}\right)
n+W_{0}-W_{1}.
\end{eqnarray*}%
This result can be also proved by mathematical induction. Note that from (a)
we get%
\begin{equation*}
\sum_{k=0}^{n}W_{2k+1}=\frac{1}{2}\left(
W_{2n+3}+W_{2n+2}+W_{2n+1}-W_{2}-W_{1}+W_{0}\right) -\sum_{k=0}^{n}W_{2k}.
\end{equation*}%
Now, (c) follows from the last equation. 
\endproof%
\end{description}

From the above Theorem we have the following Corollary which gives sum
formulas of Pell-Padovan numbers (take $W_{n}=R_{n}$ with $%
R_{0}=1,R_{1}=1,R_{2}=1$).

\begin{corollary}
For $n\geq 0,$ Pell-Padovan numbers have the following property:

\begin{description}
\item[(a)] $\sum_{k=0}^{n}R_{k}=\frac{1}{2}\left(
R_{n+3}+R_{n+2}-R_{n+1}-1\right) .$

\item[(b)] $\sum_{k=0}^{n}R_{2k}=R_{2n+1}-n.$

\item[(c)] $\sum_{k=0}^{n}R_{2k+1}=\frac{1}{2}\left(
R_{2n+3}+R_{2n+2}-R_{2n+1}+2n-1\right) .$
\end{description}
\end{corollary}

Taking $W_{n}=C_{n}$ with $C_{0}=3,C_{1}=0,C_{2}=2$ in the last Theorem, we
have the following Corollary which presents sum formulas of Pell-Perrin
numbers.

\begin{corollary}
For $n\geq 0,$ Pell-Perrin numbers have the following property:

\begin{description}
\item[(a)] $\sum_{k=0}^{n}C_{k}=\frac{1}{2}\left(
C_{n+3}+C_{n+2}-C_{n+1}+1\right) .$

\item[(b)] $\sum_{k=0}^{n}C_{2k}=C_{2n+1}-n+3.$

\item[(c)] $\sum_{k=0}^{n}C_{2k+1}=\frac{1}{2}\left(
C_{2n+3}+C_{2n+2}-C_{2n+1}+2n-5\right) .$
\end{description}
\end{corollary}

Taking $r=0,s=1,t=2$ in Theorem \ref{theorem:firstbncvg} (a) and (b) (or
(c)), we obtain the following Proposition.

\begin{proposition}
If $r=0,s=1,t=2$ then for $n\geq 0$ we have the following formulas:

\begin{description}
\item[(a)] $\sum_{k=0}^{n}W_{k}=\allowbreak \frac{1}{2}\left(
W_{n+3}+W_{n+2}-W_{2}-W_{1}\right) .$

\item[(b)] $\sum_{k=0}^{n}W_{2k}=\frac{1}{2}\left(
W_{2n+1}+2W_{2n}-W_{1}\right) .$

\item[(c)] $\sum_{k=0}^{n}W_{2k+1}=\frac{1}{2}\left(
W_{2n+2}+2W_{2n+1}-W_{2}\right) .$
\end{description}
\end{proposition}

Taking $W_{n}=Q_{n}$ with $Q_{0}=1,Q_{1}=1,Q_{2}=1$ in the last Proposition,
we have the following Corollary which presents linear sum formulas of
Jacobsthal-Padovan numbers.

\begin{corollary}
For $n\geq 0,$ Jacobsthal-Padovan numbers have the following properties.

\begin{description}
\item[(a)] $\sum_{k=0}^{n}Q_{k}=\frac{1}{2}\left( Q_{n+3}+Q_{n+2}-2\right) .$

\item[(b)] $\sum_{k=0}^{n}Q_{2k}=\frac{1}{2}\left( Q_{2n+1}+2Q_{2n}-1\right)
.$

\item[(c)] $\sum_{k=0}^{n}Q_{2k+1}=\frac{1}{2}\left(
Q_{2n+2}+2Q_{2n+1}-1\right) .$
\end{description}
\end{corollary}

From the last Proposition, we have the following Corollary which gives
linear sum formulas of Jacobsthal-Perrin numbers (take $W_{n}=D_{n}$ with $%
D_{0}=3,D_{1}=0,D_{2}=2$).

\begin{corollary}
For $n\geq 0,$ Jacobsthal-Perrin numbers have the following properties.

\begin{description}
\item[(a)] $\sum_{k=0}^{n}D_{k}=\frac{1}{2}\left( D_{n+3}+D_{n+2}-2\right) .$

\item[(b)] $\sum_{k=0}^{n}D_{2k}=\frac{1}{2}\left( D_{2n+1}+2D_{2n}\right) .$

\item[(c)] $\sum_{k=0}^{n}D_{2k+1}=\frac{1}{2}\left(
D_{2n+2}+2D_{2n+1}-2\right) .$
\end{description}
\end{corollary}

Taking $r=1,s=0,t=1$ in Theorem \ref{theorem:firstbncvg} (a) and (c), we
obtain the following Proposition.

\begin{proposition}
If $r=1,s=0,t=1$ then for $n\geq 0$ we have the following formulas:

\begin{description}
\item[(a)] $\sum_{k=0}^{n}W_{k}=W_{n+3}-W_{2}.$

\item[(b)] $\sum_{k=0}^{n}W_{2k}=\frac{1}{3}%
(W_{2n+2}+W_{2n+1}+2W_{2n}-W_{2}-W_{1}+W_{0})\allowbreak .$

\item[(c)] $\sum_{k=0}^{n}W_{2k+1}=\frac{1}{3}%
(2W_{2n+2}+2W_{2n+1}+W_{2n}-2W_{2}+W_{1}-W_{0})\allowbreak .$
\end{description}
\end{proposition}

From the last Proposition, we have the following Corollary which presents
linear sum formulas of Narayana numbers (take $W_{n}=N_{n}$ with $%
N_{0}=0,N_{1}=1,N_{2}=1$).

\begin{corollary}
For $n\geq 0,$ Narayana numbers have the following properties.

\begin{description}
\item[(a)] $\sum_{k=0}^{n}N_{k}=N_{n+3}-1.$

\item[(b)] $\sum_{k=0}^{n}N_{2k}=\frac{1}{3}(N_{2n+2}+N_{2n+1}+2N_{2n}-2).$

\item[(c)] $\sum_{k=0}^{n}N_{2k+1}=\frac{1}{3}%
(2N_{2n+2}+2N_{2n+1}+N_{2n}-1). $
\end{description}
\end{corollary}

Taking $r=1,s=1,t=2$ in Theorem \ref{theorem:firstbncvg} (a) and (c), we
obtain the following Proposition.

\begin{proposition}
If $r=1,s=1,t=2$ then for $n\geq 0$ we have the following formulas:

\begin{description}
\item[(a)] $\sum_{k=0}^{n}W_{k}=\frac{1}{3}(W_{n+3}-W_{n+1}-W_{2}+W_{0}).$

\item[(b)] $\sum_{k=0}^{n}W_{2k}=\frac{1}{3}(W_{2n+1}+2W_{2n}-W_{1}+W_{0})%
\allowbreak .$

\item[(c)] $\sum_{k=0}^{n}W_{2k+1}=\frac{1}{3}%
(W_{2n+2}+2W_{2n+1}-W_{2}+W_{1}).$
\end{description}
\end{proposition}

Taking $W_{n}=J_{n}^{(3)}$ with $J_{0}^{(3)}=0,J_{1}^{(3)}=1,J_{2}^{(3)}=1$
in the last Proposition, we have the following Corollary which presents
linear sum formulas of third order Jacobsthal numbers.

\begin{corollary}
For $n\geq 0,$ third order Jacobsthal numbers have the following properties.

\begin{description}
\item[(a)] $\sum_{k=0}^{n}J_{k}^{(3)}=\frac{1}{3}%
(J_{n+3}^{(3)}-J_{n+1}^{(3)}-1).$

\item[(b)] $\sum_{k=0}^{n}J_{2k}^{(3)}=\frac{1}{3}%
(J_{2n+1}^{(3)}+2J_{2n}^{(3)}-1).$

\item[(c)] $\sum_{k=0}^{n}J_{2k+1}^{(3)}=\frac{1}{3}%
(J_{2n+2}^{(3)}+2J_{2n+1}^{(3)}).$
\end{description}
\end{corollary}

From the last Proposition, we have the following Corollary which gives
linear sum formulas of third order Jacobsthal-Lucas numbers (take $%
W_{n}=j_{n}$ with $j_{0}^{(3)}=2,j_{1}^{(3)}=1,j_{2}^{(3)}=5$).

\begin{corollary}
For $n\geq 0,$ third order Jacobsthal-Lucas numbers have the following
properties.

\begin{description}
\item[(a)] $\sum_{k=0}^{n}j_{k}^{(3)}=\frac{1}{3}%
(j_{n+3}^{(3)}-j_{n+1}^{(3)}-3).$

\item[(b)] $\sum_{k=0}^{n}j_{2k}^{(3)}=\frac{1}{3}%
(j_{2n+1}^{(3)}+2j_{2n}^{(3)}+1).$

\item[(c)] $\sum_{k=0}^{n}j_{2k+1}^{(3)}=\frac{1}{3}%
(j_{2n+2}^{(3)}+2j_{2n+1}^{(3)}-4).$
\end{description}
\end{corollary}

\section{Sum formulas of Generalized Tribonacci Numbers with Negative
Subscripts}

The following Theorem presents some linear summing formulas (identities) of
generalized Tribonacci numbers with negative subscripts.

\begin{theorem}
\label{theorem:tyuhnbcfgxbzeqwa}For $n\geq 1,$ we have the following
formulas:

\begin{description}
\item[(a)] (Sum of the generalized Tribonacci numbers with negative indices)
If $r+s+t-1\neq 0,$ then 
\begin{equation*}
\sum_{k=1}^{n}W_{-k}=\frac{-(r+s+\allowbreak
t)W_{-n-1}-(s+t)W_{-n-2}-tW_{-n-3}+W_{2}+(1-r)W_{1}+(1-r-s)W_{0}}{r+s+t-1}.
\end{equation*}

\item[(b)] If $\allowbreak \left( r+s+t-1\right) \left( r-s+t+1\right) \neq
0 $ then%
\begin{equation*}
\sum_{k=1}^{n}W_{-2k}=\frac{%
\begin{array}{c}
-(r+t)W_{-2n+1}+(r^{2}+rt+s-1)W_{-2n}+(st-t)W_{-2n-1} \\ 
+(1-s)W_{2}+(t+rs)W_{1}+(1-rt-2s-r^{2}+s^{2})W_{0}%
\end{array}%
}{\left( r+s+t-1\right) \left( r-s+t+1\right) }
\end{equation*}%
and 
\begin{equation*}
\sum_{k=1}^{n}W_{-2k+1}=\frac{%
\begin{array}{c}
(s-1)W_{-2n+1}-(t+rs)W_{-2n}-(t^{2}+rt)W_{-2n-1} \\ 
+(r+t)W_{2}+(1-r^{2}-rt-s)W_{1}+(t-st)W_{0}%
\end{array}%
}{\left( r+s+t-1\right) \left( r-s+t+1\right) }.
\end{equation*}

\item[(c)] If $\allowbreak \left( r+s+t-1\right) \left( r-s+t+1\right) \neq
0\wedge r+t=0\wedge s\neq 1$ then 
\begin{equation*}
\sum_{k=1}^{n}W_{-2k}=\frac{-W_{-2n}-tW_{-2n-1}+W_{2}+tW_{1}+(1-s)W_{0}}{s-1}
\end{equation*}%
and%
\begin{equation*}
\sum_{k=1}^{n}W_{-2k+1}=\frac{1}{s-1}\left(
-W_{-2n+1}-tW_{-2n}+W_{1}+tW_{0}\right) .
\end{equation*}%
Note that (c) is a special case of (b).
\end{description}
\end{theorem}

\textit{Proof.} \ 

\begin{description}
\item[(a)] Using the recurrence relation%
\begin{equation*}
W_{-n+3}=rW_{-n+2}+sW_{-n+1}+tW_{-n}\Rightarrow W_{-n}=-\frac{s}{t}%
W_{-(n-1)}-\frac{r}{t}W_{-(n-2)}+\frac{1}{t}W_{-(n-3)}
\end{equation*}%
i.e.%
\begin{equation*}
tW_{-n}=W_{-n+3}-rW_{-n+2}-sW_{-n+1}
\end{equation*}%
or%
\begin{equation*}
W_{-n}=\frac{1}{t}W_{-n+3}-\frac{r}{t}W_{-n+2}-\frac{s}{t}W_{-n+1}
\end{equation*}%
we obtain%
\begin{eqnarray*}
tW_{-n} &=&W_{-n+3}-rW_{-n+2}-sW_{-n+1} \\
tW_{-n+1} &=&W_{-n+4}-rW_{-n+3}-sW_{-n+2} \\
tW_{-n+2} &=&W_{-n+5}-rW_{-n+4}-sW_{-n+3} \\
&&\vdots \\
tW_{-2} &=&W_{1}-r\times W_{0}-s\times W_{-1} \\
tW_{-1} &=&W_{2}-r\times W_{1}-s\times W_{0}.
\end{eqnarray*}%
If we add the above equations by side by, we get%
\begin{equation*}
\sum_{k=1}^{n}W_{-k}=\frac{%
\begin{array}{c}
-(rW_{-n-1}+s(W_{-n-1}+W_{-n-2})+t(W_{-n-1}+W_{-n-2}+W_{-n-3}) \\ 
-W_{2}+(r-1)W_{1}+(r+s-1)W_{0})%
\end{array}%
}{r+s+t-1}.
\end{equation*}

\item[(b) and (c)] Using the recurrence relation%
\begin{equation*}
W_{-n+3}=rW_{-n+2}+sW_{-n+1}+tW_{-n}
\end{equation*}%
i.e.%
\begin{equation*}
sW_{-n+1}=W_{-n+3}-rW_{-n+2}-tW_{-n}
\end{equation*}%
we obtain%
\begin{eqnarray*}
sW_{-2n+1} &=&W_{-2n+3}-rW_{-2n+2}-tW_{-2n} \\
sW_{-2n+3} &=&W_{-2n+5}-rW_{-2n+4}-tW_{-2n+2} \\
&&\vdots \\
sW_{-3} &=&W_{-1}-rW_{-2}-tW_{-4} \\
sW_{-1} &=&W_{1}-rW_{0}-tW_{-2}.
\end{eqnarray*}%
If we add the equations by side by, we get%
\begin{equation}
s\sum_{k=1}^{n}W_{-2k+1}=(-W_{-2n+1}+W_{1}+%
\sum_{k=1}^{n}W_{-2k+1})-r(-W_{-2n}+W_{0}+\sum_{k=1}^{n}W_{-2k})-t(%
\sum_{k=1}^{n}W_{-2k}).  \label{equati:nvrtsfcdzxcazs}
\end{equation}%
Similarly, using the recurrence relation%
\begin{equation*}
W_{-n+3}=rW_{-n+2}+sW_{-n+1}+tW_{-n}
\end{equation*}%
i.e.%
\begin{equation*}
sW_{-n+1}=W_{-n+3}-rW_{-n+2}-tW_{-n}
\end{equation*}%
we obtain%
\begin{eqnarray*}
sW_{-2n} &=&W_{-2n+2}-rW_{-2n+1}-tW_{-2n-1} \\
sW_{-2n+2} &=&W_{-2n+4}-rW_{-2n+3}-tW_{-2n+1} \\
&&\vdots \\
sW_{-6} &=&W_{-4}-rW_{-5}-tW_{-7} \\
sW_{-4} &=&W_{-2}-rW_{-3}-tW_{-5} \\
sW_{-2} &=&W_{0}-rW_{-1}-tW_{-3}.
\end{eqnarray*}%
If we add the above equations by side by, we get%
\begin{equation*}
s\sum_{k=1}^{n}W_{-2k}=(-W_{-2n}+W_{0}+\sum_{k=1}^{n}W_{-2k})-r(%
\sum_{k=1}^{n}W_{-2k+1})-t(W_{-2n-1}-W_{-1}+\sum_{k=1}^{n}W_{-2k+1}).
\end{equation*}%
Since%
\begin{equation*}
W_{-1}=(-\frac{s}{t}W_{0}-\frac{r}{t}W_{1}+\frac{1}{t}W_{2}).
\end{equation*}%
it follows that%
\begin{eqnarray}
s\sum_{k=1}^{n}W_{-2k}
&=&(-W_{-2n}+W_{0}+\sum_{k=1}^{n}W_{-2k})-r(\sum_{k=1}^{n}W_{-2k+1})
\label{equati:mbvuosdczdxa} \\
&&-t(W_{-2n-1}-(-\frac{s}{t}W_{0}-\frac{r}{t}W_{1}+\frac{1}{t}%
W_{2})+\sum_{k=1}^{n}W_{-2k+1}).  \notag
\end{eqnarray}%
Then, solving system (\ref{equati:nvrtsfcdzxcazs})-(\ref{equati:mbvuosdczdxa}%
) the required results of (b) and (c) follow. 
\endproof%
\end{description}

Note that (c) of the above theorem can be written as follows: If $%
r+t=0\wedge s\neq 1$ then 
\begin{equation*}
\sum_{k=1}^{n}W_{-2k}=\frac{-W_{-2n}+rW_{-2n-1}+W_{2}-rW_{1}+(1-s)W_{0}}{s-1}
\end{equation*}%
and%
\begin{equation*}
\sum_{k=1}^{n}W_{-2k}=\frac{-W_{-2n}+rW_{-2n-1}+W_{2}-rW_{1}+(1-s)W_{0}}{s-1}%
.
\end{equation*}

Next, we present several sum formulas (identities).

Taking $r=s=t=1$ in Theorem \ref{theorem:tyuhnbcfgxbzeqwa} (a) and (b), we
obtain the following Proposition.

\begin{proposition}
If $r=s=t=1$ then for $n\geq 1$ we have the following formulas:

\begin{description}
\item[(a)] $\sum_{k=1}^{n}W_{-k}=\allowbreak \frac{1}{2}\left(
-3W_{-n-1}-2W_{-n-2}-W_{-n-3}+W_{2}-W_{0}\right) .$

\item[(b)] $\sum_{k=1}^{n}W_{-2k}=\frac{1}{2}\left(
-W_{-2n+1}+W_{-2n}+W_{1}-W_{0}\right) .$

\item[(c)] $\sum_{k=1}^{n}W_{-2k+1}=\allowbreak \frac{1}{2}\left(
-W_{-2n}-W_{-2n-1}+W_{2}-W_{1}\right) .$
\end{description}
\end{proposition}

From the above Proposition, we have the following Corollary which gives
linear sum formulas of Tribonacci numbers (take $W_{n}=T_{n}$ with $%
T_{0}=0,T_{1}=1,T_{2}=1$).

\begin{corollary}
For $n\geq 1,$ Tribonacci numbers have the following properties.

\begin{description}
\item[(a)] $\sum_{k=1}^{n}T_{-k}=\allowbreak \frac{1}{2}%
(-3T_{-n-1}-2T_{-n-2}-T_{-n-3}+1).$

\item[(b)] $\sum_{k=1}^{n}T_{-2k}=\allowbreak \frac{1}{2}%
(-T_{-2n+1}+T_{-2n}+1).$

\item[(c)] $\sum_{k=1}^{n}T_{-2k+1}=\allowbreak \frac{1}{2}%
(-T_{-2n}-T_{-2n-1}).$
\end{description}
\end{corollary}

Taking $W_{n}=K_{n}$ with $K_{0}=3,K_{1}=1,K_{2}=3$ in the above
Proposition, we have the following Corollary which gives linear sum formulas
of Tribonacci-Lucas numbers.

\begin{corollary}
For $n\geq 1,$ Tribonacci-Lucas numbers have the following properties:

\begin{description}
\item[(a)] $\sum_{k=1}^{n}K_{-k}=\allowbreak \frac{1}{2}%
(-3K_{-n-1}-2K_{-n-2}-K_{-n-3}).$

\item[(b)] $\sum_{k=1}^{n}K_{-2k}=\allowbreak \frac{1}{2}%
(-K_{-2n+1}+K_{-2n}-2).$

\item[(c)] $\sum_{k=1}^{n}K_{-2k+1}=\allowbreak \frac{1}{2}%
(-K_{-2n}-K_{-2n-1}+2).$
\end{description}
\end{corollary}

Taking $r=2,s=1,t=1$ in Theorem \ref{theorem:tyuhnbcfgxbzeqwa} (a) and (b),
we obtain the following Proposition.

\begin{proposition}
If $r=2,s=1,t=1$ then for $n\geq 1$ we have the following formulas:

\begin{description}
\item[(a)] $\sum_{k=1}^{n}W_{-k}=\frac{1}{3}\left(
-4W_{-n-1}-2W_{-n-2}-W_{-n-3}+W_{2}-W_{1}-2W_{0}\right) .$

\item[(b)] $\sum_{k=1}^{n}W_{-2k}=\allowbreak \frac{1}{3}\left(
-W_{-2n+1}+2W_{-2n}+W_{1}-2W_{0}\right) .$

\item[(c)] $\sum_{k=1}^{n}W_{-2k+1}=\allowbreak \frac{1}{3}\left(
-W_{-2n}-W_{-2n-1}+W_{2}-2W_{1}\right) .$
\end{description}
\end{proposition}

From the last Proposition, we have the following Corollary which gives
linear sum formulas of third-order Pell numbers (take $W_{n}=P_{n}^{(3)}$
with $P_{0}^{(3)}=0,P_{1}^{(3)}=1,P_{2}^{(3)}=2$).

\begin{corollary}
For $n\geq 1,$ third-order Pell numbers have the following properties.

\begin{description}
\item[(a)] $\sum_{k=1}^{n}P_{-k}^{(3)}=\frac{1}{3}%
(-4P_{-n-1}^{(3)}-2P_{-n-2}^{(3)}-P_{-n-3}^{(3)}+1).$

\item[(b)] $\sum_{k=1}^{n}P_{-2k}^{(3)}=\frac{1}{3}%
(-P_{-2n+1}^{(3)}+2P_{-2n}^{(3)}+1).$

\item[(c)] $\sum_{k=1}^{n}P_{-2k+1}^{(3)}=\allowbreak \frac{1}{3}%
(-P_{-2n}^{(3)}-P_{-2n-1}^{(3)}).$
\end{description}
\end{corollary}

Taking $W_{n}=Q_{n}^{(3)}$ with $Q_{0}^{(3)}=3,Q_{1}^{(3)}=2,Q_{2}^{(3)}=6$
in the last Proposition, we have the following Corollary which gives linear
sum formulas of third-order Pell-Lucas numbers.

\begin{corollary}
For $n\geq 1,$ third-order Pell-Lucas numbers have the following properties.

\begin{description}
\item[(a)] $\sum_{k=1}^{n}Q_{-k}^{(3)}=\frac{1}{3}%
(-4Q_{-n-1}^{(3)}-2Q_{-n-2}^{(3)}-Q_{-n-3}^{(3)}-2).$

\item[(b)] $\sum_{k=1}^{n}Q_{-2k}^{(3)}=\frac{1}{3}%
(-Q_{-2n+1}^{(3)}+2Q_{-2n}^{(3)}-4).$

\item[(c)] $\sum_{k=1}^{n}Q_{-2k+1}^{(3)}=\allowbreak \frac{1}{3}%
(-Q_{-2n}^{(3)}-Q_{-2n-1}^{(3)}+2).$
\end{description}
\end{corollary}

From the last Proposition, we have the following Corollary which presents
linear sum formulas of third-order modified Pell numbers (take $%
W_{n}=E_{n}^{(3)}$ with $E_{0}^{(3)}=0,E_{1}^{(3)}=1,E_{2}^{(3)}=1$).

\begin{corollary}
For $n\geq 1,$ third-order modified Pell numbers have the following
properties.

\begin{description}
\item[(a)] $\sum_{k=1}^{n}E_{-k}^{(3)}=\frac{1}{3}%
(-4E_{-n-1}^{(3)}-2E_{-n-2}^{(3)}-E_{-n-3}^{(3)}).$

\item[(b)] $\sum_{k=1}^{n}E_{-2k}^{(3)}=\frac{1}{3}%
(-E_{-2n+1}^{(3)}+2E_{-2n}^{(3)}+1).$

\item[(c)] $\sum_{k=1}^{n}E_{-2k+1}^{(3)}=\allowbreak \frac{1}{3}%
(-E_{-2n}^{(3)}-E_{-2n-1}^{(3)}-1).$
\end{description}
\end{corollary}

Taking $r=0,s=1,t=1$ in Theorem \ref{theorem:tyuhnbcfgxbzeqwa} (a) and (b),
we obtain the following Proposition.

\begin{proposition}
If $r=0,s=1,t=1$ then for $n\geq 1$ we have the following formulas:

\begin{description}
\item[(a)] $\sum_{k=1}^{n}W_{-k}=-2W_{-n-1}-2W_{-n-2}-W_{-n-3}+W_{2}+W_{1}.$

\item[(b)] $\sum_{k=1}^{n}W_{-2k}=-W_{-2n+1}+W_{1}.$

\item[(c)] $\sum_{k=1}^{n}W_{-2k+1}=-W_{-2n}-W_{-2n-1}+W_{2}.$
\end{description}
\end{proposition}

Taking $W_{n}=P_{n}$ with $P_{0}=1,P_{1}=1,P_{2}=1$ in the last Proposition,
we have the following Corollary which gives linear sum formulas of Padovan
numbers.

\begin{corollary}
For $n\geq 1,$ Padovan numbers have the following properties.

\begin{description}
\item[(a)] $\sum_{k=1}^{n}P_{-k}=-2P_{-n-1}-2P_{-n-2}-P_{-n-3}+2.$

\item[(b)] $\sum_{k=1}^{n}P_{-2k}=-P_{-2n+1}+1.$

\item[(c)] $\sum_{k=1}^{n}P_{-2k+1}=-P_{-2n}-P_{-2n-1}+1.$
\end{description}
\end{corollary}

From the last Proposition, we have the following Corollary which presents
linear sum formulas of Perrin numbers (take $W_{n}=E_{n}$ with $%
E_{0}=3,E=0,E_{2}=2$).

\begin{corollary}
For $n\geq 1,$ Perrin numbers have the following properties.

\begin{description}
\item[(a)] $\sum_{k=1}^{n}E_{-k}=-2E_{-n-1}-2E_{-n-2}-E_{-n-3}+2.$

\item[(b)] $\sum_{k=1}^{n}E_{-2k}=-E_{-2n+1}.$

\item[(c)] $\sum_{k=1}^{n}E_{-2k+1}=-E_{-2n}-E_{-2n-1}+2.$
\end{description}
\end{corollary}

Taking $W_{n}=S_{n}$ with $S_{0}=0,S_{1}=0,S_{2}=1$ in the last Proposition,
we have the following Corollary which gives linear sum formulas of
Padovan-Perrin numbers.

\begin{corollary}
For $n\geq 1,$ Padovan-Perrin numbers have the following properties.

\begin{description}
\item[(a)] $\sum_{k=1}^{n}S_{-k}=-2S_{-n-1}-2S_{-n-2}-S_{-n-3}+1.$

\item[(b)] $\sum_{k=1}^{n}S_{-2k}=-S_{-2n+1}.$

\item[(c)] $\sum_{k=1}^{n}S_{-2k+1}=-S_{-2n}-S_{-2n-1}+1.$
\end{description}
\end{corollary}

If $r=0,s=2,t=1$ then $\left( r+s+t-1\right) \left( r-s+t+1\right)
=\allowbreak 0$ so we can't use Theorem \ref{theorem:tyuhnbcfgxbzeqwa} (b)
and (c). In other words, the method of the proof Theorem \ref%
{theorem:tyuhnbcfgxbzeqwa} (b) and (c) can't be used to find $%
\sum_{k=0}^{n}W_{2k}$ and $\sum_{k=0}^{n}W_{2k+1}.$ Therefore we need
another method to find them which is given in the following Theorem.

\begin{theorem}
If $r=0,s=2,t=1$ then for $n\geq 1$ we have the following formulas:

\begin{description}
\item[(a)] $\sum_{k=1}^{n}W_{-k}=\frac{1}{2}\left(
-3W_{-n-1}-3W_{-n-2}-W_{-n-3}+W_{2}+W_{1}-W_{0}\right) .$

\item[(b)] $%
\sum_{k=1}^{n}W_{-2k}=-W_{-2n+1}+W_{-2n}+(W_{1}-W_{0})+(W_{2}-W_{1}-W_{0})n.$

\item[(c)] $\sum_{k=1}^{n}W_{-2k+1}=\frac{1}{2}%
(W_{-2n+1}-3W_{-2n}-W_{-2n-1}+(W_{2}-W_{1}+W_{0})+2(-W_{2}+W_{1}+W_{0})n).$
\end{description}
\end{theorem}

Proof.

\begin{description}
\item[(a)] Taking $r=0,s=2,t=1$ in Theorem \ref{theorem:tyuhnbcfgxbzeqwa}
(a) we obtain (a).

\item[(b) and (c)] Proof can be done as in the proof of Theorem \ref%
{theorem:bayrambhospd}. Induction also can be used for the proof.
\end{description}

From the last Theorem, we have the following Corollary which gives sum
formula of Pell-Padovan numbers (take $W_{n}=R_{n}$ with $%
R_{0}=1,R=1,R_{2}=1 $).

\begin{corollary}
For $n\geq 1,$ Pell-Padovan numbers have the following property:

\begin{description}
\item[(a)] $\sum_{k=1}^{n}R_{-k}=\frac{1}{2}\left(
-3R_{-n-1}-3R_{-n-2}-R_{-n-3}+1\right) .$

\item[(b)] $\sum_{k=1}^{n}R_{-2k}=-R_{-2n+1}+R_{-2n}-n.$

\item[(c)] $\sum_{k=1}^{n}R_{-2k+1}=\frac{1}{2}%
(R_{-2n+1}-3R_{-2n}-R_{-2n-1}+1+2n).$
\end{description}
\end{corollary}

Taking $W_{n}=C_{n}$ with $C_{0}=3,C=0,C_{2}=2$ in the last Theorem, we have
the following Corollary which gives sum formulas of Pell-Perrin numbers.

\begin{corollary}
For $n\geq 1,$ Pell-Perrin numbers have the following property:

\begin{description}
\item[(a)] $\sum_{k=1}^{n}C_{-k}=\frac{1}{2}\left(
-3C_{-n-1}-3C_{-n-2}-C_{-n-3}-1\right) $

\item[(b)] $\sum_{k=1}^{n}C_{-2k}=-C_{-2n+1}+C_{-2n}-3-n$

\item[(c)] $\sum_{k=1}^{n}C_{-2k+1}=\frac{1}{2}%
(C_{-2n+1}-3C_{-2n}-C_{-2n-1}+5+2n)$
\end{description}
\end{corollary}

Taking $r=0,s=1,t=2$ in Theorem \ref{theorem:tyuhnbcfgxbzeqwa} (a) and (b),
we obtain the following Proposition.

\begin{proposition}
If $r=0,s=1,t=2$ then for $n\geq 1$ we have the following formulas:

\begin{description}
\item[(a)] $\sum_{k=1}^{n}W_{-k}=\frac{1}{2}\left(
-3W_{-n-1}-3W_{-n-2}-2W_{-n-3}+W_{2}+W_{1}\right) .$

\item[(b)] $\sum_{k=1}^{n}W_{-2k}=\frac{1}{2}\left( -W_{-2n+1}+W_{1}\right)
. $

\item[(c)] $\sum_{k=1}^{n}W_{-2k+1}=\frac{1}{2}\left(
-W_{-2n}-2W_{-2n-1}+W_{2}\right) .$
\end{description}
\end{proposition}

From the last Proposition, we have the following Corollary which gives
linear sum formulas of Jacobsthal-Padovan numbers (take $W_{n}=Q_{n}$ with $%
Q_{0}=1,Q_{1}=1,Q_{2}=1$).

\begin{corollary}
For $n\geq 1,$ Jacobsthal-Padovan numbers have the following properties.

\begin{description}
\item[(a)] $\sum_{k=1}^{n}Q_{-k}=\frac{1}{2}\left(
-3Q_{-n-1}-3Q_{-n-2}-2Q_{-n-3}+2\right) .$

\item[(b)] $\sum_{k=1}^{n}Q_{-2k}=\frac{1}{2}\left( -Q_{-2n+1}+1\right) .$

\item[(c)] $\sum_{k=1}^{n}Q_{-2k+1}=\frac{1}{2}\left(
-Q_{-2n}-2Q_{-2n-1}+1\right) .$
\end{description}
\end{corollary}

Taking $W_{n}=D_{n}$ with $D_{0}=3,D_{1}=0,D_{2}=2$ in the last Proposition,
we have the following Corollary which gives linear sum formulas of
Jacobsthal-Perrin numbers.

\begin{corollary}
For $n\geq 1,$ Jacobsthal-Perrin numbers have the following properties.

\begin{description}
\item[(a)] $\sum_{k=1}^{n}D_{-k}=\frac{1}{2}\left(
-3D_{-n-1}-3D_{-n-2}-2D_{-n-3}+2\right) .$

\item[(b)] $\sum_{k=1}^{n}D_{-2k}=\frac{-1}{2}D_{-2n+1}.$

\item[(c)] $\sum_{k=1}^{n}D_{-2k+1}=\frac{1}{2}\left(
-D_{-2n}-2D_{-2n-1}+2\right) .$
\end{description}
\end{corollary}

Taking $r=1,s=0,t=1$ in Theorem \ref{theorem:tyuhnbcfgxbzeqwa}, we obtain
the following Proposition.

\begin{proposition}
If $r=1,s=0,t=1$ then for $n\geq 1$ we have the following formulas:

\begin{description}
\item[(a)] $\sum_{k=1}^{n}W_{-k}=-2W_{-n-1}-W_{-n-2}-W_{-n-3}+W_{2}.$

\item[(b)] $\sum_{k=1}^{n}W_{-2k}=\frac{1}{3}\left(
-2W_{-2n+1}+W_{-2n}-W_{-2n-1}+W_{2}+W_{1}-W_{0}\right) .$

\item[(c)] $\sum_{k=1}^{n}W_{-2k+1}=\frac{1}{3}\left(
-W_{-2n+1}-W_{-2n}-2W_{-2n-1}+2W_{2}-W_{1}+W_{0}\right) .$
\end{description}
\end{proposition}

From the above Proposition, we have the following Corollary which gives
linear sum formulas of Narayana numbers (take $W_{n}=N_{n}$ with $%
N_{0}=0,N_{1}=1,N_{2}=1).$

\begin{corollary}
For $n\geq 1,$ Narayana numbers have the following properties.

\begin{description}
\item[(a)] $\sum_{k=1}^{n}N_{-k}=-2N_{-n-1}-N_{-n-2}-N_{-n-3}+1.$

\item[(b)] $\sum_{k=1}^{n}N_{-2k}=\frac{1}{3}\left(
-2N_{-2n+1}+N_{-2n}-N_{-2n-1}+2\right) .$

\item[(c)] $\sum_{k=1}^{n}N_{-2k+1}=\frac{1}{3}\left(
-N_{-2n+1}-N_{-2n}-2N_{-2n-1}+1\right) \allowbreak .$
\end{description}
\end{corollary}

Taking $r=1,s=1,t=2$ in Theorem \ref{theorem:tyuhnbcfgxbzeqwa}, we obtain
the following Proposition.

\begin{proposition}
If $r=1,s=1,t=2$ then for $n\geq 1$ we have the following formulas:

\begin{description}
\item[(a)] $\sum_{k=1}^{n}W_{-k}=\frac{1}{3}%
(-4W_{-n-1}-3W_{-n-2}-2W_{-n-3}+W_{2}-W_{0}).$

\item[(b)] $\sum_{k=1}^{n}W_{-2k}=\frac{1}{3}\left(
-W_{-2n+1}+W_{-2n}+W_{1}-W_{0}\right) .$

\item[(c)] $\sum_{k=1}^{n}W_{-2k+1}=\frac{1}{3}\left(
-W_{-2n}-2W_{-2n-1}+W_{2}-W_{1}\right) .$
\end{description}
\end{proposition}

Taking $W_{n}=J_{n}$ with $J_{0}=0,J_{1}=1,J_{2}=1$ in the last Proposition,
we have the following Corollary which gives linear sum formulas of third
order Jacobsthal numbers.

\begin{corollary}
For $n\geq 1,$ third order Jacobsthal numbers have the following properties.

\begin{description}
\item[(a)] $\sum_{k=1}^{n}J_{-k}^{(3)}=\frac{1}{3}%
(-4J_{-n-1}^{(3)}-3J_{-n-2}^{(3)}-2J_{-n-3}^{(3)}+1).$

\item[(b)] $\sum_{k=1}^{n}J_{-2k}^{(3)}=\frac{1}{3}%
(-J_{-2n+1}^{(3)}+J_{-2n}^{(3)}+1).$

\item[(c)] $\sum_{k=1}^{n}J_{-2k+1}^{(3)}=\frac{1}{3}%
(-J_{-2n}^{(3)}-2J_{-2n-1}^{(3)})\allowbreak .$
\end{description}
\end{corollary}

From the last Proposition, we have the following Corollary which gives
linear sum formulas of third order Jacobsthal-Lucas numbers (take $%
W_{n}=j_{n}^{(3)}$ with $j_{0}^{(3)}=2,j_{1}^{(3)}=1,j_{2}^{(3)}=5$).

\begin{corollary}
For $n\geq 1,$ third order Jacobsthal-Lucas numbers have the following
properties.

\begin{description}
\item[(a)] $\sum_{k=1}^{n}j_{-k}^{(3)}=\frac{1}{3}%
(-4j_{-n-1}^{(3)}-3j_{-n-2}^{(3)}-2j_{-n-3}^{(3)}+3).$

\item[(b)] $\sum_{k=1}^{n}j_{-2k}^{(3)}=\frac{1}{3}%
(-j_{-2n+1}^{(3)}+j_{-2n}^{(3)}-1).$

\item[(c)] $\sum_{k=1}^{n}j_{-2k+1}^{(3)}=\frac{1}{3}%
(-j_{-2n}^{(3)}-2j_{-2n-1}^{(3)}+4)\allowbreak .$
\end{description}
\end{corollary}

\end{document}